\documentclass[11pt]{article}
\usepackage[backrefs, alphabetic,initials]{amsrefs}
\usepackage{amsfonts,amsmath,amscd,amssymb,amsthm}
\usepackage{amsfonts,amsmath,amscd,amssymb,amsthm,latexcad}
\usepackage{fullpage}
\usepackage{mathrsfs}
\def\0{{\bar 0}}
\def\ww{\wedge}

\def\1{{\bar 1}}

\newcommand{\CQ}{\overline{\cQ}}

\def\Ext{{\operatorname{Ext}}}

\newcommand{\rra}{\longleftarrow}
\newcommand{\dra}{\Rightarrow}
\newcommand{\dla}{\Leftarrow}
\newcommand{\itema}{\item[{{\rm(a)}}]}
\newcommand{\itemb}{\item[{{\rm(b)}}]}
\newcommand{\itemc}{\item[{{\rm(c)}}]}

\newcommand{\gL}{\Lambda}

\newcommand{\noi}{\noindent}
\newcommand{\ga}{\alpha}

\def\Supp{{\operatorname{Supp}}}

\newcommand{\gC}{\Gamma}

\newcommand{\ff}{\footnote}
\newfont{\eufm}{eufm10 scaled\magstep1}

\newcommand{\cD}{\mathcal{D}}

\newcommand{\cA}{\mathcal{A}}

\newcommand{\cP}{\mathcal{P}}

\newcommand{\cL}{\mathcal{L}}

\newcommand{\cQ}{\mathcal{Q}}

\newcommand{\bco}{\begin{conjecture}}
\newcommand{\ba}{\begin{alg}}
\newcommand{\ea}{\end{alg}}
\newcommand{\eco}{\end{conjecture}}
\newcommand{\bpf}{\begin{proof}}
\newcommand{\epf}{\end{proof}}
\newcommand{\bt}{\begin{theorem}}
\newcommand{\et}{\end{theorem}}

\newcommand{\br}{\begin{rem}}
\newcommand{\er}{\end{rem}}
\newcommand{\brs}{\begin{rems}}
\newcommand{\ers}{\end{rems}}
\newcommand{\bi}{\begin{itemize}}
\newcommand{\bl}{\begin{lemma}}
\newcommand{\bsul}{\begin{sublemma}}
\newcommand{\esul}{\end{sublemma}}
\newcommand{\bp}{\begin{proposition}}
\newcommand{\be}{\begin{equation}}
\newcommand{\bc}{\begin{corollary}}
\newcommand{\bexa}{\begin{example}}
\newcommand{\eexa}{\end{example}}
\newcommand{\bex}{\begin{exercise}}
\newcommand{\eex}{\end{exercise}}
\newcommand{\btab}{\begin{tab}}
\newcommand{\etab}{\end{tab}}
\newcommand{\ei}{\end{itemize}}
\newcommand{\el}{\end{lemma}}
\newcommand{\ep}{\end{proposition}}
\newcommand{\ee}{\end{equation}}
\newcommand{\ec}{\end{corollary}}
\newcommand{\Bc}{\begin{center}}
\newcommand{\Ec}{\end{center}}
\newcommand{\bh}{\begin{hyp}}
\newcommand{\eh}{\end{hyp}}
\newcommand{\bhs}{\begin{hyps}}
\newcommand{\ehs}{\end{hyps}}

\numberwithin{equation}{section}%
\newcommand{\lra}{\longrightarrow}

\begin{document}
\title{Table of Contents}
\newcommand{\rh}{\epsilon^{\widehat{\rho}}}
\newcommand{\pii}{\prod_{n=1}^\infty}
\newcommand{\sii}{\sum_{n= - \infty}^\infty}
\newtheorem{theorem}{Theorem}[section]
 \newtheorem{lemma}[theorem]{Lemma}
 \newtheorem{example}[theorem]{Example}
  \newtheorem{sublemma}[theorem]{Sublemma}
   \newtheorem{corollary}[theorem]{Corollary}
 \newtheorem{conjecture}[theorem]{Conjecture}
 \newtheorem{hyps}[theorem]{Hypotheses}
  \newtheorem{hyp}[theorem]{Hypothesis}
 \newtheorem{alg}[theorem]{Algorithm}
  \newtheorem{sico}[theorem]{Sign Convention}
 \newtheorem{tab}[theorem]{\quad \quad \quad \quad \quad \quad \quad \quad \quad \quad \quad \quad \quad Table }
  \newtheorem{rem}[theorem]{Remark}
  \newtheorem{rems}[theorem]{Remarks}
 \newtheorem{proposition}[theorem]{Proposition}
\font\twelveeufm=eufm10 scaled\magstep1 \font\teneufm=eufm10
\font\nineeufm=eufm9 \font\eighteufm=eufm8 \font\seveneufm=eufm7
\font\sixeufm=eufm6 \font\fiveeufm=eufm5
\newfam\eufmfam
\textfont\eufmfam=\twelveeufm \scriptfont\eufmfam=\nineeufm
\scriptscriptfont\eufmfam=\sixeufm \textfont\eufmfam=\teneufm
\scriptfont\eufmfam=\seveneufm \scriptscriptfont\eufmfam\fiveeufm
\newtheorem{exercise}{}
\setcounter{exercise}{0} \numberwithin{exercise}{section}
\newenvironment{emphit}{\begin{theorem}}{\end{theorem}}
\setcounter{Thm}{0} \numberwithin{Thm}{section} 
%

\title{
The lattice of submodules of a multiplicity free module.}
\author{Ian M. Musson\ff{Research supported by  NSA Grant H98230-12-1-0249.} \\Department of Mathematical Sciences\\
University of Wisconsin-Milwaukee\\ email: {\tt
musson@uwm.edu}}
\maketitle

\begin{abstract}
In this paper we determine, under some mild restrictions, the lattice of submodules $\gL$ of a module $M$ all of whose composition factors have multiplicity one.
Such a lattice is distributive, and hence determined by its poset of down-sets $P$. We define a directed Ext graph $\Ext_\gL$ of  $\gL$ and show that if  $\Ext_\gL$ is acyclic, then $\Ext_\gL$  determines $P$.
The result applies to multiplicity free indecomposable modules for finite dimensional algebras with acyclic Ext graph. It also applies to some deformed Verma modules which arise in the Jantzen sum formula basic classical simple Lie superalgebras in the deformed case. \end{abstract}
\section{Introduction.}
Two basic problems in representation theory are to determine the simple objects in some abelian category  $\cA$, and then to determine the nontrivial extensions between simples.
Often that is as much as we can expect to say about the internal structure of objects in $\cA$, even those of finite length. For example  it is possible for a Verma module to contain infinitely many submodules \cite{CD}. However if we consider multiplicity free objects then we can often say more.  For convenience we will assume that $\cA$ is a category of (finite length) modules over some ring.\\ \\
In some situations the lattice of submodules $\gL$ of a module $M$ can be determined completely.  Indeed if $M$ is multiplicity free, then $\gL$ is distributive, so is determined by the Fundamental Theorem on Distributive Lattices in terms of its poset $P$ of join irreducible submodules, see Theorem \ref{dog}.
 We that if  the underlying graph of $\Ext_\gL$ is acyclic, then $\Ext_\gL$  determines $P$. \\ \\
 In Section 3 we give some applications of our result.
Recall that a {\it quiver} is a 4-tuple ${\cQ}
= (\cQ_0, \cQ_1, s, t)$ where $\cQ_0, \cQ_1$ are finite sets of {\it
vertices} and {\it arrows} respectively and $s,t: \cQ_1
\longrightarrow \cQ_0$ are maps assigning to each arrow its {\it
source} and {\it target} respectively. We call the graph whose vertex and edge set are $\cQ_0, \cQ_1$ as the {\it underlying graph $\CQ$ of $\cQ$}. If its underlying graph has no cycles, we say that $\CQ$ is acyclic. 
In section 2 we refer to  quivers as  digraphs since this is more common in the world of combinatorics.

The result on lattices applies to indecomposable multiplicity free module over path algebras $K\cQ$ provided that $\CQ$ is acyclic.
In addition the result applies to certain deformed Verma modules for a
classical simple Lie superalgebra. In these cases $\gL$ is the free distributive lattice $\cD_k$ freely generated by $k$ join irreducibles \cite{M2}. \\ \\
I am grateful to Peter Cameron for many pleasant conversations about lattices.  Another approach to the lattice of submodules of a multiplicity free module was given by Alperin \cite{A}.  I thank Zongzhu Lin for pointing this out.

\section{Background from  lattice theory} \label{ali}
A finite poset $\gL$ is  a {\it lattice} if every pair $A, B$ of elements of $\gL$ has a greatest lower bound $A \ww B$ and a least upper bound $A\vee B$. These are necessarily unique. Furthermore the conditions $B\le A,$ $B=A \ww B$ and $A=A\vee B$ are equivalent. In this case we say that $[A,B]$ is an {\it interval}. If $B< A$ and $B\le C\le A$ implies that $C=A$ or $C= B$ we say that  $A$ {\it covers} $B$. 
We say that $A$ is {\it join irreducible} if $A =X\vee Y$ implies that either $A=X$ or $A= Y,$ and
 We say that the interval $[X,Y]$ is {\it uniserial} if there is a unique maximal chain from $X$ to $Y$. Let $\equiv$ be the smallest equivalence relation on the set of intervals such that $[A\vee B,B] \equiv [A,A\ww B]$. The equivalence classes under $\equiv$ will be called {\it simple lattice factors} of $\gL$. We say that $A$ is {\it join irreducible} if $A =X\vee Y$ implies that either $A=X$ or $A= Y.$
We say that $\gL$ is {\it multiplicity free} if given a maximal chain
\be \label{nut} 0=X_0 < X_1 < \ldots < X_p = 1 \ee
in $\gL,$ $[X_{i-1},X_i]\equiv [X_{j-1},X_j]$ implies that $i=j.$ In general we use upper case letters to denote elements of $\gL$ and lower case letters for simple lattice factors. We sometimes refer to elements of $\gL$ as {\it submodules}.
\\ \\
We will apply our result to the case where $\gL$ is the lattice of submodules of a module $M$ of finite length. So we assume henceforth that $\gL$ is  modular. If $M$ is multiplicity free as a module then $\gL$ is multiplicity free. The converse of this statement is false, as shown by the cyclic group of order 4.\\ \\
The following is well known, but we include a proof for completeness.
\bl A multiplicity free modular lattice is distributive \el
\bpf
By \cite{B} Theorem II.13 it is enough to show that $\gL$ contains no sublattice with Hasse diagram

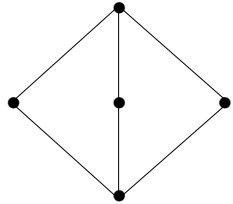
\begin{picture}(104,100)(-135,50)
\thinlines
\put(42,64.0){$\bullet$}
\put(42,135.0){$\bullet$}
\put(42,99.0){$\bullet$}
\put(2,99.0){$\bullet$}
\put(82,99.0){$\bullet$}
\drawpath{44.0}{138.0}{4.0}{102.0}
\drawpath{86.0}{102.0}{45.0}{138.0}
\drawpath{4.0}{102.0}{44.0}{66.0}
\drawpath{86.0}{102.0}{45.0}{66.0}
\drawpath{44.5}{67.0}{44.5}{137.0}
\end{picture}

\noi This follows since $\gL$ is  multiplicity free.\epf
\noi Note that every element of $\gL$ is a join of irreducible elements.
In addition we have
\bl \label{cat}\bi \itema If $A$ is join irreducible element then the join of all elements of $\gL$ that are strictly less than $A$ is the unique element that is covered by $A$.
\itemb  If $X$ covers $Y$, then $[X,Y] \equiv [A,B]$ where $A$ is join irreducible.
\itemc  If $[X,Y]$ is uniserial, then $[X,Y] \equiv [A,B]$ where $A$ is join irreducible.
\ei
\el \bpf (a) is immediate, and (b) is a special case of (c), so we prove (c).
Suppose that $A$ is minimal such $[X,Y] \equiv [A,B]$ for some $B$. We claim that $A$ is join irreducible. If not then $A = U\vee V$ where $U<A$ and $V<A$.
If also $U\vee B = A$ we have
$[A,B] = [U\vee B,B]\equiv [U,U\ww B]$ and this contradicts the minimality of $A$. Hence $U'=U\vee B$, and similarly $V'=V\vee B = B$ are both greater or equal to $B$ and strictly less than $A$. Since
$A = U'\vee V',$ and $[A,B]$ is uniserial, this is impossible.
\epf
\noi If $U$ covers $V$ and $X$ covers $Y$ we write $[U,V]\dra [X,Y]$ or $[X,Y] \dla [U,V]$ if $U \vee Y = X,$ and $U \ww Y = V.$
\bl \label{owl}\bi \itema The relation
$\dra$ is transitive.
\itemb If $X$ is join irreducible and $[U,V]\dra [X,Y]$, then $[U,V]=[X,Y]$.\ei\el \bpf
Suppose  $[U,V]\dra [X,Y]$ as above, and  $[A,B] \dra [U,V]$.
Then $U = A\vee V$ and $B = A\ww V.$ Thus $A \vee Y= A\vee V\vee Y=U\vee Y = X,$ and $A \ww Y \le A\ww U\ww Y=  A\ww V = B,$ so $[A,B]\dra [X,Y]$.
\\ \\
To prove (b)
Suppose $U \vee Y = X,$ and $U \ww Y = V.$ Since $X$ is join irreducible, and $X \neq Y$ we have $U = X,$ and $V= X \ww Y = Y.$
\epf
\bl \label{pig}If $\gL$ is distributive, the join irreducible in Lemma \ref{cat} (b) is unique.\el
\bpf If not, then by Lemma \ref{pig} (b) we can find a sequence of intervals
$$[A,B], [X_0,Y_0],
[A_1,B_1],\ldots [A_m,B_m],[X_m,Y_m]$$ such that
$$[A,B]\dra [X_0,Y_0]\dla[A_1,B_1]\dra [X_1,Y_1]\dla\ldots \dla [A_m,B_m]\dra [X_m,Y_m]\dla [X,Y]$$ with
$A$ and $X$ join irreducible. Choose such a sequence with $m$ minimal, and set $U=A\ww A_1, V=B\ww B_1.$
Since $\gL$ is distributive, it follows  that $[A,B]\dla [U,V]\dra[A_1,B_1].$ Hence by  Lemma \ref{pig}
$[A,B]= [U,V]$, so by transitivity  $[A,B]\dra [X_1,Y_1]$.  We have now found a sequence of shorter length, a contradiction.
\epf
We will be a bit  lazy and write $[U,V]=y$ to mean that $U$ covers $V$, and the equivalence class of $[U,V]$ is $y$. 
If $X$ is join irreducible, we denote the unique maximal submodule and simple factor module of $X$ by $X^0$ and $x$ respectively. If there is any chance of ambiguity we will say $X$ is the join irreducible with top $x$.
\bl \label{on} Any submodule of $\gL$ with $y$ as a simple lattice  factor contains $Y$.\el
\bpf If $X$ contains $y$ as a simple lattice  factor, choose a submodule of $X$ having $y$ as a simple lattice  factor such that no proper submodule has $y$ as a simple lattice  factor.  This submodule is join irreducible, so  by uniqueness it is equal to $Y$.
\epf

\noi If $P$ is a poset then a subset $I$ of $P$ is called an {\it down set} or {\it order ideal}  if whenever  $y\le x$ and $x\in I$ we have $y \in I.$ If $I$ and $J$ are down sets, then so are $I\cap J$ and $I\cup J$. Thus the set of down sets of $P$ forms a distributive lattice $J(P)$ with $\ww,\vee$ given by $\cap, \cup$  respectively.\\ \\
The next result is known as the fundamental theorem of distributive lattices.
\bt \label{dog} Let $\gL$ be a finite distributive lattice, and let
$P$ be the subposet of $\gL$ consisting of join irreducibles.  Then there is a unique $($up to isomorphism$)$ finite poset $P$ such that $\gL \cong J(P)$.\et \bpf  This is shown in \cite{S} Theorem 3.4.1.\epf
\noi Note that from any poset $P$, we obtain an acyclic digraph $\cP$, with vertex set $P$ and an edge $x\lra y$ iff $x$ covers $y.$
Moreover any acyclic digraph can be obtained in this way. We call $\cP$ the {\it digraph associated to} $P$.\\ \\
\noi
Now suppose that $x,y$ are simple lattice factors.  If there is a length two interval $[u,w]$ containing a unique proper subinterval $[v,w]$  such that $[u,v] \equiv x$ and $[v,w] \equiv y,$ we say that $[u,w]$ is a non-trivial extension of $y$ by $x$, and write $\Ext_{\gL}(x,y)\neq 0.$ 
\bl \label{cow}Suppose $X,Y$ are join irreducible with unique simple factors $x,y$ respectively. Let  $P$ be the  poset of join irreducible ideals of $\gL$.
\bi \itema If $X$ covers $Y$  in $P$, then  $\Ext_{\gL}(x,y)\neq 0$.
\itemb Suppose $\Ext_{\gL}(x,y)\neq 0$, and the edge $x\lra y$ is the only path from $x$ to $y$ in the
the underlying graph of the directed $\Ext_\gL$ graph.
 Then $X$ covers $Y$  in $P.$
\itemc Suppose
the underlying graph of the directed $\Ext$ graph of $\gL$ is acylic. Then $X$ covers $Y$  in $P$ iff $\Ext_{\gL}(x,y)\neq 0$
\ei
\el
\bpf To prove (a) suppose $X$ covers $Y$, that is $X, Y$ are join irreducibles such that there is no join irreducible strictly between $X$ and $Y$. Choose a submodule $A$ of $X^o$ that is maximal such that $A \cap Y = Y^o$.  We claim that $A+Y =X^o$. Obviously $A+Y \subseteq X^o$.  If the inclusion is proper
then in the interval $[X^o,A+Y]$ we can write $X^o$ as a join of join irreducibles. This
contradicts the assumption that $X$ covers $Y.$ Now set
$C = X/A$.  Then $C$ has a submodule $X^o/A \cong y$ with factor module isomorphic to $x$. It follows that $X/A$ is uniserial and  $\Ext_{\gL}(x,y)\neq 0$.
\\ \\
Since (c) follows at once from (a) and (b), it remains to prove (b). Suppose $\Ext_{\gL}(x,y)\neq 0$. By Lemmas \ref{cat} and \ref{pig} there is a uniserial interval  $[X,U]$ in $\gL$ which is a non-trivial extension of $y$ by $x$. 
Suppose $X$ does not cover $Y$ in the poset of join irreducibles, and
 let   $X > Z_1>\ldots >Z_n >Y$ be a  maximal chain with $n\ge 1.$ Then  by (a) there are arrows $x\lra z_1\lra \ldots \lra z_n \lra y.$ However this gives a cycle in the Ext graph.  \epf
\bt \label{tom}Let $P$ be the poset of join irreducible elements of the distributive lattice $\gL$, and let $\cP$ be the digraph with vertex set the set of simple lattice factors of $\gL$
and with an edge $x\lra y$ iff
$\Ext_{\gL}(x,y)\neq 0$. If  $\cP$ is acyclic, then $\cP$ is isomorphic to the associated digraph of $P$. \et
\bpf By Lemmas
\ref{pig} and \ref{cow} the map sending $X$ to $x$ is an isomorphism.\epf

\section{Applications and Examples.}
Let $M$ be a multiplicity free  $R$-module of finite length, and let $\gL$ be its lattice of submodules. In this case the composition factors of $M$ are the same as the simple lattice factors of $\gL$.
By definition the $\Ext_{\gL}$ graph (resp. the $\Ext_{R}$ graph) of $M$ has the composition factors of
$M$ as its vertices, with an edge from $x$ to $y$ iff 
 $\Ext_{\gL}(x,y) \neq 0$ (resp.  $\Ext_{R}(x,y) \neq 0$).
From the point of view of ring theory $\Ext_{R}$  is more natural than
$\Ext_{\gL},$ so we are interested in conditions where they are the same.
\bl \label{hen}If $M$ is an indecomposable module, then  $\Ext_\gL$ is connected. \el
\bpf Suppose
$\Ext_{\gL}$ is a disjoint union
$\Ext_\gL = \gC_1 \cup \gC_2$, and for $i=1,2$ let $M_i$ be the largest submodule of $M$ all of whose composition factors are vertices of $\gC_i$.  Then $M= M_1\oplus M_2$. \epf
\subsection{Finite Dimensional Algebras}
 A {\it representation} of the quiver $\cQ$ is a family of vector spaces $M_x$ indexed by $x
\in \cQ_0$ together with maps
\[ f_\ga : M_{s(\ga)} \longrightarrow M_{t(\ga)} \]
for every arrow $\ga \in \cQ_1$.  Given such a representation, $M =\bigoplus_{x\in \cQ} M_x$ is a
$K\cQ$-module and every $K\cQ$-module arises in this way.
The dimension vector $\dim M$ of $M$
is the vector with entries indexed by $\cQ_0$ such that the entry
corresponding to $x$ is $\dim M_x$. Clearly  $M$ is multiplicity free iff every entry of $\dim M$ is equal to zero or one. We set $\Supp \;M = \{x\in \cQ_0|M_x\neq 0\}.$ If $H$ is a subset of the vertex set of a digraph $G$, the {\it subgraph of $G$ induced by} $H$ has vertex set $H$ and has an arrow from $x$ to $y$ iff there is an arrow  from $x$ to $y$ in $G.$
\bt Suppose that $\cQ$ is a finite acyclic quiver and set $R=K\cQ$. If $M$ is an indecomposable  multiplicity free module, with lattice of submodules $\gL$,
then $\Ext_\gL$  is the subgraph of  $\Ext_R$ induced by $\Supp M.$   \et
\bpf By \cite{ASS} Lemma II.2.5, $R$ is indecomposable as an algebra iff $\cQ$ is connected.  Now if $e$ is any central idempotent in $R$ we have $M = eM \oplus (1-e)M.$  Thus we can assume that $\cQ$ is connected, that is $\cQ$ is a tree.
Thus if any edge of $\cQ$ is removed (without removing any vertices) the result is a disconnected graph.
Suppose that $x$ and $y$ are composition factors of $M$ and $\Ext_R(x,y)\neq 0.$ Then there is an edge $\ga$ from $x$ to $y$ in $\cQ$. Since $M$ is indecomposable this implies that the map $f_\ga:M_x\lra M_y$ is also non-zero.  Thus $\Ext_\gL(x,y)\neq 0.$ \epf
\bc Let $K$ be an algebraically closed field, and $R$ a finite dimensional $K$-algebra. Suppose that the Ext quiver of $R$ is acyclic.  If $M$ is an indecomposable  multiplicity free module, with lattice of submodules $\gL$,
then $\Ext_\gL$  is the subgraph of  $\Ext_R$ induced by $\Supp M.$   \ec \bpf     
 If $B$ is the basic algebra associated to $R$  \cite{ASS} I.6, then $R$ is Morita equivalent to $B$, so we can replace $R$ by $B$ to assume at the outset that $R$ is basic. Let $\cQ$ be the directed quiver of $R$. Then by  \cite{ASS} Theorem II.3.7 and Lemma III.2.12,
there is an admissible ideal $I$ of the path algebra $K\cQ$ such that $R \cong K\cQ/I.$ Thus $M$ is also a $K\cQ$-module, so the result follows. \epf
\brs {\rm \bi \itema
There are situations, for example modular group algebras, where an indecomposable projective $M$  has isomorphic socle and cosocle.  However it is clear that the  to determine the lattice of submodules, it is enough to consider $rad(M)/soc(M)$ instead.
    \itemb The following simple example may be instructive.  Let $K$ be a field of characteristic three and $G$ the symmetric group of degree three.  Then $e = \frac{1+(12)}{2}$  is idempotent in $R =KG$, and $P=Ae$ is uniserial with socle and cosocle isomorphic to
        the trivial module $V^+$, and has  $rad(P)/soc(P)$ isomorphic to
        the sign module $V^-.$ The Ext graph of $R$ is
        \[V^+ \stackrel{\textstyle{\lra}}{\rra} V^-.\]
Therefore if $M=rad(P)$ and $\gL$ is the lattice of submodules of $M$, then $\Ext_{A}(V^-,V^+)\neq 0$ but $\Ext_{\gL}(V^-,V^+)=0$.
\ei} \ers
\subsection{Deformed Verma modules for basic classical Lie superalgebras.}

Let $\{P_1,\ldots , P_n\}$ be a set of $n$ independent propositions, and $\cD = \cD_n$ the set of propositions that can be formed from the $P_i$ using conjunctions and disjunctions.  We put a partial order on $\cD$ by declaring that $P\le \cQ$ means that $P\vee \cQ =P.$ We include 0 (resp.) 1 as elements of $\cD_n$ which are less (resp. greater) than any other element. We call the lattice $\cD_n$ the {\it Dedekind lattice of order} $n$. The cardinality $ |\cD_n|$ is known as the $n^{th}$ Dedekind number.\\ \\
 According to Wikipedia
``the Dedekind numbers are a rapidly growing sequence of integers named after Richard Dedekind, who defined them in 1897. The Dedekind number $M(n)$ counts the number of monotonic Boolean functions of n variables. Equivalently, it counts the number of antichains of subsets of an $n$-element set, the number of elements in a free distributive lattice with $n$ generators, or the number of abstract simplicial complexes with $n$ elements."\\ \\
In the notation of Theorem \ref{dog} we have $\cD_n \cong J(P)$ for $P$ the boolean lattice with $n$ atoms. Additionally in  \cite{Do} the Dedekind lattice (without 0 and 1) is called the lattice of lattice polynomials $\cL_n$ on $n$ variables. For the Hasse diagram of $\cL_3$ see \cite{Do} Figure 6.8.1. This lattice is also isomorphic to the lattice of square free monomial ideals in $n$ variables, see \cite{MiS}. In spite of these many incarnations of the Dedekind lattice,  the value of $|\cD_n|$ is known  only for $n\le 8.$

In \cite{M2} we show that if $M$ is a deformed Verma module for a
basic classical simple Lie superalgebra, and $M$ has degree of atypicality $n$, then the lattice of submodules of $M$   is isomorphic to $\cD_n$. See \cite{M} for background on Lie superalgebras. It was this result which prompted our investigations into lattices.


\end{document}